\documentclass[letterpaper, 10 pt, conference]{ieeeconf}  
\usepackage{amsmath,amssymb}
\usepackage[pdftex]{graphicx}
\usepackage{color}
\newtheorem{thm}{\bf Theorem}
\newtheorem{prb}{\bf Problem}
\newtheorem{lem}{\bf Lemma}

\IEEEoverridecommandlockouts                              

\title{\LARGE \bf
Data-driven Estimation of the Algebraic Riccati Equation for \\the Discrete-Time Inverse Linear Quadratic Regulator Problem
}

\author{Shuhei Sugiura$^{1}$, Ryo Ariizumi$^{2}$, Masaya Tanemura$^{3}$, Toru Asai$^{1}$, and Shun-ichi Azuma$^{4}$
\thanks{*This work was supported in part by the Japan Society for the Promotion of Science KAKENHI under Grant JP22K04027 and in part by JST FOREST Program JPMJFR2123.}
\thanks{$^{1}$Shuhei Sugiura and Toru Asai are with the Graduate School of Engineering, Nagoya University, Nagoya, Japan {\tt\small sugiura.shuhei.j6@s.mail.nagoya-u.ac.jp asai@nuem.nagoya-u.ac.jp}} %
\thanks{$^{2}$Ryo Ariizumi is with the Department of Mechanical System Engineering, Tokyo University of Agriculture and Technology, Tokyo, Japan {\tt\small ryoariizumi@go.tuat.ac.jp}}%
\thanks{$^{3}$Masaya Tanemura is with the Graduate School of Science and Technology, Shinshu University, Nagano, Japan {\tt\small tanemura@shinshu-u.ac.jp}}%
\thanks{$^{2}$Shun-ichi Azuma is with the Graduate School of Informatics, Kyoto University, Kyoto, Japan {\tt\small azuma.shunichi.3e@kyoto-u.ac.jp}}%
}

\begin{document}
\onecolumn This work has been submitted to the IEEE for possible publication. Copyright may be transferred without notice, after which this version may no longer be accessible.
\twocolumn

\newpage
\maketitle
\thispagestyle{empty}
\pagestyle{empty}

\begin{abstract}
In this paper, we propose a method for estimating the algebraic Riccati equation (ARE) with respect to an unknown discrete-time system from the system state and input observation. The inverse optimal control (IOC) problem asks, ``What objective function is optimized by a given control system?'' The inverse linear quadratic regulator (ILQR) problem is an IOC problem that assumes a linear system and quadratic objective function. The ILQR problem can be solved by solving a linear matrix inequality that contains the ARE. However, the system model is required to obtain the ARE, and it is often unknown in fields in which the IOC problem occurs, for example, biological system analysis. Our method directly estimates the ARE from the observation data without identifying the system. This feature enables us to economize the observation data using prior information about the objective function. We provide a data condition that is sufficient for our method to estimate the ARE. We conducted a numerical experiment to demonstrate that our method can estimate the ARE with less data than system identification if the prior information is sufficient.
\end{abstract}
\section{INTRODUCTION}
The inverse optimal control (IOC) problem is a framework used to determine the objective function optimized by a given control system. By optimizing the objective function obtained by IOC, a different system can imitate the given control system's behavior. For example, several studies have been conducted on robots imitating human or animal movements \cite{Mombaur,Li,El-Hussieny}. To apply IOC, in these studies, the researchers assumed that human and animal movements are optimal for unknown criteria. Such an optimality assumption holds in some cases \cite{Alexander}.

The first IOC problem, which assumes a single-input linear time-invariant system and quadratic objective function, was proposed by Kalman \cite{Kalman}. Anderson \cite{Anderson} generalized the IOC problem in \cite{Kalman} to the multi-input case. For this type of IOC problem, called the inverse linear quadratic regulator (ILQR) problem, Molinari \cite{Molinari} proposed a necessary and sufficient condition for a linear feedback input to optimize some objective functions. Moylan et al. \cite{Moylan} proposed and solved the IOC problem for nonlinear systems.  IOC in \cite{Kalman,Anderson,Molinari,Moylan} is based on control theory, whereas Ng et al. \cite{Ng} proposed inverse reinforcement learning (IRL), which is IOC based on machine learning. Recently, IRL has become an important IOC framework, along with control theory methods \cite{Ab Azar}.

We consider the ILQR problem with respect to an unknown discrete-time system. Such a problem often occurs in biological system analysis, which is the main application target of IOC. The control theory approach solves the ILQR problem by solving a linear matrix inequality (LMI) that contains the algebraic Riccati equation (ARE) \cite{Priess}. However, to obtain the ARE, the system model is required. The IRL approach also has difficulty solving our problem. There is IRL with unknown systems \cite{Herman}, and IRL with continuous states and input spaces \cite{Aghasadeghi}; however, to the best of our knowledge, no IRL exists with both.

In the continuous-time case, we proposed a method to estimate the ARE from the observation data of the system state and input directly \cite{continuous}. In the present paper, we estimate the ARE for a discrete-time system using a method that extends the result in \cite{continuous}. Similarly to \cite{continuous}, our method transforms the ARE by multiplying the observed state and input on both sides. However, this technique alone cannot calculate ARE without the system model because the form of the ARE differs between continuous and discrete-time. We solve this problem using inputs from the system's controller. Also, the use of such inputs enable us to estimate the ARE without knowing the system's control gain. We prove that the equation obtained from this transformation is equivalent to the ARE if the dataset is appropriate. The advantage of our method is the economization of the observation data using prior information about the objective function. We conducted a numerical experiment to demonstrate that our method can estimate the ARE with less data than system identification if the prior information is sufficient.

The structure of the remainder of this paper is as follows: In Section II, we formulate the problem considered in this paper. In Section III, we propose our estimation method and prove that the estimated equation is equivalent to the ARE. In Section IV, we describe numerical experiments that confirm our statement in Section III. Section V concludes the paper.
\section{PROBLEM FORMULATION}
We consider the following discrete-time linear system:
\begin{equation}
\label{transition}
x{\left(k+1\right)}=Ax{\left(k\right)}+Bu{\left(k\right)},
\end{equation}
where $A\in{\mathbb R}^{n\times n}$ and $B\in{\mathbb R}^{n\times m}$ are constant matrices, and $x{\left(k\right)}\in{\mathbb R}^n$ and $u{\left(k\right)}\in{\mathbb R}^m$ are the state and the input at time $k\in\mathbb Z$, respectively.

Let $\mathbb Z_+=\left\{0,1,2,\cdots\right\}$. We write the set of real-valued symmetric $n\times n$ matrices as ${\mathbb R}^{n\times n}_{\rm sym}$. In the LQR problem \cite{book1}, we determine the input sequence $u{\left(k\right)}$ $\left(k\in\mathbb Z_+\right)$ that minimizes the following objective function:
\begin{equation}
\label{objective}
J{\left(u\right)}=\sum_{k=0}^\infty \left[x{\left(k\right)}^\top Qx{\left(k\right)}+u{\left(k\right)}^\top R u{\left(k\right)}\right],
\end{equation}
where $Q\in{\mathbb R}^{n\times n}_{\rm sym}$ and $R\in{\mathbb R}^{m\times m}_{\rm sym}$.

We assume that the system (\ref{transition}) is controllable and matrices $Q$ and $R$ are positive definite. Then, for any $x{\left(0\right)}$, the LQR problem has a unique solution written as the following linear state feedback law:
\begin{equation}
\label{231025140638}
u{\left(k\right)}=-\left(R+B^\top PB\right)^{-1}B^\top PAx{\left(k\right)},
\end{equation}
where $P\in{\mathbb R}^{n\times n}_{\rm sym}$ is a unique positive definite solution of the ARE:
\begin{equation}
\label{ARE1}
A^\top PA-P+Q-A^\top PB\left(R+B^\top PB\right)^{-1}B^\top PA=0.
\end{equation}
For $u$ defined in (\ref{231025140638}), $J{\left(u\right)}=x{\left(0\right)}^\top Px{\left(0\right)}$ holds.

In the ILQR problem, we find positive definite matrices $Q\in{\mathbb R}^{n\times n}_{\rm sym}$ and $R\in{\mathbb R}^{m\times m}_{\rm sym}$ such that an input $u{\left(k\right)}=-Kx{\left(k\right)}$ with the given gain $K\in{\mathbb R}^{m\times n}$ is a solution of the LQR problem, that is, $u$ minimizes (\ref{objective}).

We can solve the ILQR problem by solving an LMI. Let the system (\ref{transition}) be controllable. Then, $u{\left(k\right)}=-Kx{\left(k\right)}$ minimizes (\ref{objective}) if and only if a positive definite matrix $P\in{\mathbb R}^{n\times n}_{\rm sym}$ exists that satisfies (\ref{ARE1}) and the following:
\begin{equation}
\label{gain}
K=\left(R+B^\top PB\right)^{-1}B^\top PA.
\end{equation}
By transforming (\ref{ARE1}) and (\ref{gain}), we obtain the following pair of equations equivalent to (\ref{ARE1}) and (\ref{gain}):
\begin{equation}
\label{ARE2}
\begin{split}
A^\top PA-P+Q-K^\top\left(R+B^\top PB\right)K=&0,\\
B^\top PA-\left(R+B^\top PB\right)K=&0.
\end{split}
\end{equation}
Hence, we can solve the ILQR problem by determining $P$, $Q\in{\mathbb R}^{n\times n}_{\rm sym}$, and $R\in{\mathbb R}^{m\times m}_{\rm sym}$ that satisfy the following LMI:
\begin{equation}
P>0,\;Q>0,\;R>0\qquad {\rm s.t. \quad Equation\;(6)}.
\end{equation}

In biological system analysis or reverse engineering, the system model and control gain is often unknown, and the ARE (\ref{ARE2}) is not readily available. Hence, we consider the following problem of determining a linear equation equivalent to (\ref{ARE2}) using the system state and input observation:
\begin{prb}[ARE estimation problem]
\it Consider controllable system (\ref{transition}) and controller $u{\left(k\right)}=-Kx{\left(k\right)}$ with unknown $A$, $B$, and $K$. Suppose $N_d$ observation data $\left(x_i{\left(0\right)},u_i,x_i{\left(1\right)}\right)$ $\left(i\in\left\{1,\ldots,N_d\right\}\right)$ of the system state and input are given, where
\begin{equation}
\label{transition2}
x_i{\left(1\right)}=Ax_i{\left(0\right)}+Bu_i\quad\left(i\in\left\{1,\ldots,N_d\right\}\right).
\end{equation}
Let $N'_d\le N_d$. $N'_d$ inputs in the data are obtained from the unknown controller as follows:
\begin{equation}
\label{asm_prb1}
u_i=-Kx_i{\left(0\right)}\quad\left(i\in\left\{1,\ldots,N'_d\right\}\right).
\end{equation}
\color{black} Determine a linear equation of $P$, $Q\in {\mathbb R}^{n\times n}_{\rm sym}$, and $R\in {\mathbb R}^{m\times m}_{\rm sym}$ with the same solution space as (\ref{ARE2}).
\end{prb}
We discuss without assuming that the data is sequential, that is, always satisfies $x_{i+1}{\left(0\right)}=x_i{\left(1\right)}$. However, in an experiment, we show that our result can also be applied to one sequential data. The ARE has multiple solutions and we call the set of these solutions, which is a linear subspace, a solution space.
\section{DATA-DRIVEN ESTIMATION OF ARE}
The simplest solution to Problem 1 is to identify the control system, that is, matrices $A$, $B$ and $K$. We define matrices $X\in\mathbb R^{n\times N_d}$, $X'\in\mathbb R^{n\times N'_d}$,  $U\in\mathbb R^{m\times N_d}$, $U'\in\mathbb R^{m\times N'_d}$, and $D\in\mathbb R^{n+m\times N_d}$ as follows:
\begin{equation}
\label{XUD}
\begin{split}
X{\left(k\right)}=&\left[x_1{\left(k\right)}\;\cdots\;x_{N_d}{\left(k\right)}\right],\\
X'{\left(k\right)}=&\left[x_1{\left(k\right)}\;\cdots\;x_{N'_d}{\left(k\right)}\right]\quad\left(k\in\left\{0,1\right\}\right),\\
U=&\left[u_1\;\cdots\;u_{N_d}\right],\quad U'=\left[u_1\;\cdots\;u_{N'_d}\right],\\
D=&\left[X{\left(0\right)}^\top\;U^\top\right]^\top
\end{split}
\end{equation}
Let the matrices $D$ and $X'{\left(0\right)}$ have row full rank. Then, matrices $A$, $B$, and $K$ are identified using the least square method as follows:
\begin{equation}
\label{SI}
\begin{split}
\left[A\;B\right]=&X{\left(1\right)}D^\top\left[DD^\top\right]^{-1}\\
K=&-U'X'{\left(0\right)}^\top\left[X'{\left(0\right)}X'{\left(0\right)}^\top\right]^{-1}\\
\end{split}
\end{equation}

\color{black}In ILQR problems, prior information about matrices $Q$ and $R$ may be obtained. However, such information cannot be used for system identification. We propose a novel method that uses the prior information and estimates the ARE using less observation data than system identification.

The following theorem provides an estimation of the ARE:
\begin{thm}
\it Consider Problem 1. Let $F$ be an $N'_d\times N_d$ matrix whose $i$th row $j$th column element $f_{i,j}$ is defined as
\begin{equation}
\label{fij}
f_{i,j}=x_i{\left(1\right)}^\top Px_j{\left(1\right)}+x_i{\left(0\right)}^\top \left(Q-P\right)x_j{\left(0\right)}+u_i^\top Ru_j.
\end{equation}
Then, the following condition is necessary for the ARE (\ref{ARE2}) to hold.
\begin{equation}
\label{estimation}
F=0.
\end{equation}
\end{thm}
\begin{proof}
We define $G_1\in\mathbb R^{n\times n}$ and $G_2\in\mathbb R^{m\times n}$ as
\begin{equation}
\label{G1G2}
\begin{split}
G_1=&A^\top PA-P+Q-K^\top\left(R+B^\top PB\right)K,\\
G_2=&B^\top PA-\left(R+B^\top PB\right)K.
\end{split}
\end{equation}
The ARE (\ref{ARE2}) is equivalent to the combination of $G_1=0$ and $G_2=0$. Let $i\in\left\{1,\ldots,N'_d\right\}$ and $j\in\left\{1,\ldots,N_d\right\}$ be arbitrary natural numbers. From the ARE, we obtain
\begin{equation}
\label{230831161420}
\begin{split}
\left[\begin{array}{*{20}{c}}
x_i{\left(0\right)}\\
{u_i}
\end{array}\right]^\top
\left[\begin{array}{*{20}{c}}
G_1&G_2^\top\\
G_2&0
\end{array}\right]
\left[\begin{array}{*{20}{c}}
x_j{\left(0\right)}\\
{u_j}
\end{array}\right]=0
\end{split}
\end{equation}
Using the assumption (\ref{transition2}) and (\ref{asm_prb1}), we can transform the left-hand side of (\ref{230831161420}) into the following:
\begin{equation}
\label{230831164105}
\begin{split}
&\left[\begin{array}{*{20}{c}}
x_i{\left(0\right)}\\
{u_i}
\end{array}\right]^\top
\left[\begin{array}{*{20}{c}}
G_1&G_2^\top\\
G_2&0
\end{array}\right]
\left[\begin{array}{*{20}{c}}
x_j{\left(0\right)}\\
{u_j}
\end{array}\right]\\
=&x_i{\left(0\right)}^\top\left(A^\top PA-P+Q-K^\top\left(R+B^\top PB\right)K\right)x_j{\left(0\right)}\\
&+x_i{\left(0\right)}^\top\left(A^\top PB-K^\top\left(R+B^\top PB\right)\right)u_j\\
&+u_i^\top\left(B^\top PA-\left(R+B^\top PB\right)K\right)x_j{\left(0\right)}\\
=&\left(Ax_i{\left(0\right)}+Bu_i\right)^\top P\left(Ax_j{\left(0\right)}+Bu_j\right)-u_i^\top B^\top PBu_j\\
&-\left(Kx_i{\left(0\right)}+u_i\right)^\top\left(R+B^\top PB\right)\left(Kx_j{\left(0\right)}+u_j\right)\\
&+u_i^\top\left(R+B^\top PB\right)u_j+x_i{\left(0\right)}^\top \left(Q-P\right)x_j{\left(0\right)}\\
=&f_{i,j}.
\end{split}
\end{equation}
Therefore, we obtain $f_{i,j}=0$ from (\ref{230831161420}) and (\ref{230831164105}), which proves Theorem 1.
\end{proof}

Equation (\ref{estimation}) is the estimated linear equation that we propose in this paper and it can be obtained directly from the observation data without system identification. Our method obtains one scalar linear equation of $P$, $Q$, and $R$ from a pair of observation data. However, at least one of the paired data must be a transition by a linear feedback input with the given gain $K$. Equation (\ref{estimation}) consists of $N_dN'_d$ scalar equations obtained in such a manner. Therefore, it is expected that if $N_d$ and $N'_d$ are sufficiently large, (\ref{estimation}) is equivalent to (\ref{ARE2}), which is also linear for the matrices $P$, $Q$, and $R$.

Equation (\ref{fij}) is similar to the Bellman equation \cite{book2}. Suppose that the ARE (\ref{ARE2}) holds, that is, $u{\left(k\right)}=-Kx{\left(k\right)}$ minimizes $J{\left(u\right)}$. Then, the following Bellman equation holds:
\begin{equation}
\label{bellman}
\begin{split}
x{\left(k\right)}^\top Px{\left(k\right)}=&x{\left(k+1\right)}^\top Px{\left(k+1\right)}+x{\left(k\right)}^\top Qx{\left(k\right)}\\
&+x{\left(k\right)}^\top K^\top RKx{\left(k\right)}.
\end{split}
\end{equation}
The equation $f_{i,j}=0$ is equivalent to the Bellman equation if $x_i{\left(0\right)}=x_j{\left(0\right)}$ and $u_i=u_j=-Kx_i{\left(0\right)}$. The novelty of Theorem 1 is that $f_{i,j}=0$ holds under the weaker condition, that is, only $u_i=-Kx_i{\left(0\right)}$. Theorem 1 is also proved from the Bellman equation (\ref{bellman}). See the supplemental material for the details.

The following theorem demonstrates that if the data satisfies the condition required for system identification, our estimation is equivalent to the ARE:
\begin{thm}
\it Consider Problem 1. Suppose that matrices $D$ and $X'{\left(0\right)}$ defined in (\ref{XUD}) have row full rank. Then, the estimation (\ref{estimation}) is equivalent to the ARE (\ref{ARE2}).
\end{thm}
\begin{proof}
From Theorem 1, the ARE (\ref{ARE2}) implies the estimation (\ref{estimation}), that is, $F=0$. As shown in the proof of Theorem 1, matrix $F$ is expressed as follows using the matrices defined in (\ref{XUD}) and (\ref{G1G2}):
\begin{equation}
\label{240202142322}
F=\left[\begin{array}{*{20}{c}}
X'{\left(0\right)}\\
U'
\end{array}\right]^\top
\left[\begin{array}{*{20}{c}}
G_1&G_2^\top\\
G_2&0
\end{array}\right]D,
\end{equation}
Because the matrix $D$ has row full rank, $F=0$ implies
\begin{equation}
\label{230920104930}
\begin{split}
X'{\left(0\right)}^\top G_1+U'^\top G_2&=0,\\
X'{\left(0\right)}^\top G_2^\top&=0.
\end{split}
\end{equation}
Because $X'{\left(0\right)}$ has row full rank, we obtain $G_1=0$ and $G_2=0$ from (\ref{230920104930}), which proves Theorem 2.
\end{proof}

The advantage of our method is the data economization that results from having prior information about $Q$ and $R$. Suppose that some elements of $Q$ and $R$ are known to be zero in advance, and independent ones are fewer than the scalar equations in the ARE (\ref{ARE2}). Then, there is a possibility that our method can estimate the ARE with fewer data than $m+n$ required for system identification.

For fixed $N_d$, our method can provide the largest number of scalar equations when $N'_d=\min\left\{n,N_d\right\}$. We use the following lemma to explain the reason:
\begin{lem}
\it Suppose that there exist $k\in\left\{1,\ldots,N'_d\right\}$ and $c_i\in\mathbb R$ $\left(i\in\left\{1,\ldots,N'_d\right\}\backslash\left\{k\right\}\right)$ that satisfy
\begin{equation}
d_k=\sum_{i=1,i\ne k}^{N'_d}c_id_i,
\end{equation}
where $d_i=\left[x_i{\left(0\right)}^\top\;u_i^\top\right]^\top$ for any $i\in\left\{1,\ldots,N_d\right\}$. Then, the following holds for any $j\in\left\{1,\ldots,N_d\right\}$:
\begin{equation}
f_{k,j}=\begin{cases}\displaystyle\sum_{i=1,i\ne k}^{N'_d}c_if_{i,j}&\left(j\ne k\right)\\
\displaystyle\sum_{i=1,i\ne k}^{N'_d}\sum_{j=1,j\ne k}^{N'_d}c_ic_jf_{i,j}&\left(j=k\right)
\end{cases}.
\end{equation}
\end{lem}
\begin{proof}
We can readily prove Lemma 1 using following derived from (\ref{240202142322})
\begin{equation}
\begin{split}
&f_{k,j}=d_k^\top
\left[\begin{array}{*{20}{c}}
G_1&G_2^\top\\
G_2&0
\end{array}\right]d_j\\
&\left(k\in\left\{1,\ldots,N'_d\right\},j\in\left\{1,\ldots,N_d\right\}\right).
\end{split}
\end{equation}
\end{proof}
Lemma 1 means that if data $d_k$ is the linear combination of other data, the equations obtained from $d_k$ are also linear combinations of other equations not using $d_k$. Hence, without noise, such a data as $d_k$ is meaningless. If $N'_d>n$, at least one of $d_1,\ldots,d_{N'_d}$ is the linear combination of other data because the inputs in those data are obtained from the same gain $K$. Therefore, $N'_d$ larger than $n$ reduces data efficiency. From the above and $f_{i,j}=f_{j,i}$, our method can provide $\frac{n\left(n+1\right)}{2}+n\left(N_d-n\right)$ equations if $N_d\ge n$.

For an example of prior information, we consider diagonal $Q$ and $R$. In this case, the number of the independent elements of $P$, $Q$, and $R$ is $\frac{n\left(n+1\right)}{2}+n+m$. Hence, if $N_d$ is $N_{d\rm min}{\left(n,m\right)}=n+1+\lceil\frac{m}{n}\rceil$ or larger, we have more equations than the independent elements. \color{black}Fig. 1 compares the numbers $N_{d\rm min}$ and $n+m$ of data required for our method and system identification, respectively, if $Q$ and $R$ are diagonal. From Fig. 1, our method may estimate the ARE with less data than system identification. Additionally, this tendency becomes stronger as the number $m$ of inputs becomes larger than the number $n$ of states.

\begin{figure}[h]
\centering
\includegraphics[width=7.5cm]{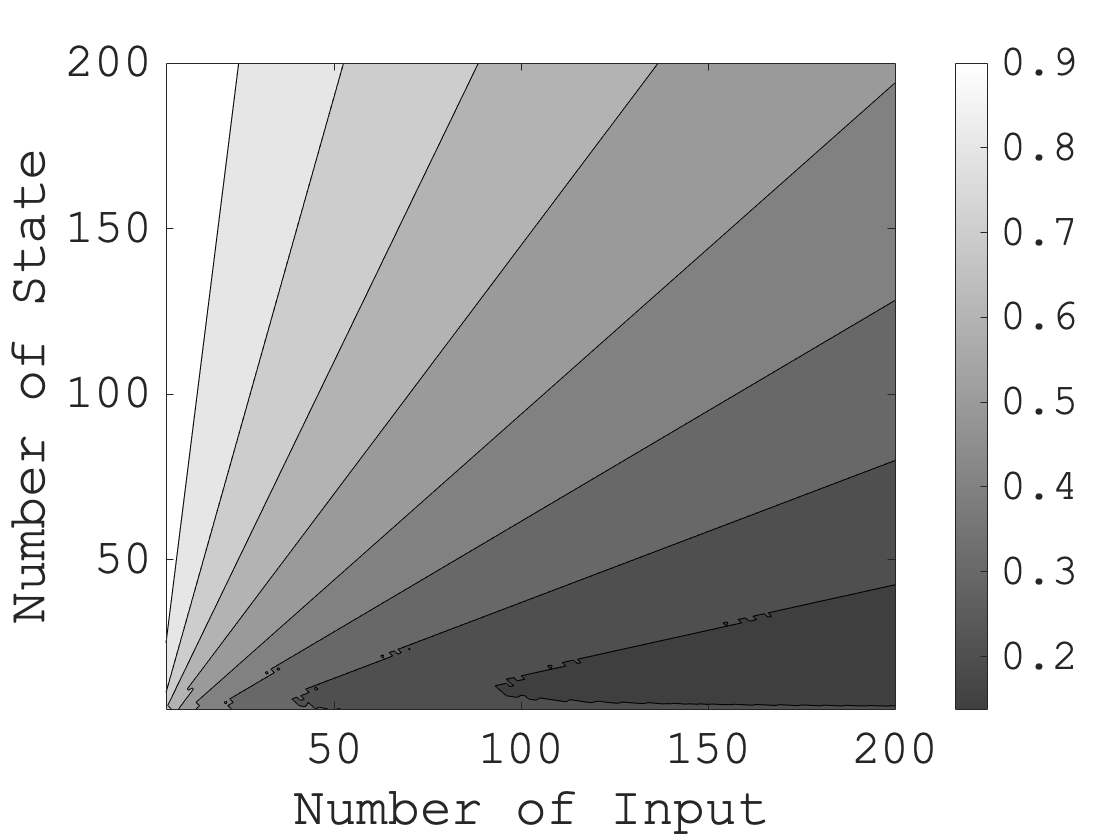}
\caption{Ratio $\frac{N_{d\rm min}}{n+m}$ when $5\le n\le 200$ and $5\le m\le 200$. The contour lines are not smooth, but this is not a mistake.}
\end{figure}
In the case of noisy data, we discuss the unbiasedness of the coefficients in (\ref{estimation}). The coefficient of $k$th row $l$th column element $q_{k,l}$ of $Q$ in (\ref{fij}) is expressed as follows:
\begin{equation}
\label{coefficientQ}
\begin{cases}x_{i,k}{\left(0\right)}x_{j,l}{\left(0\right)}+x_{i,l}{\left(0\right)}x_{j,k}{\left(0\right)}&\left(i\ne j\right)\\
x_{i,k}{\left(0\right)}x_{i,l}{\left(0\right)}&\left(i=j\right),
\end{cases}
\end{equation}
where $x_{i,k}{\left(0\right)}$ is the $k$th element of $x_i{\left(0\right)}$. Suppose that the data is a sum of the true value and zero-mean noise distributed independently for each data. Then, we can readily confirm that if $i\ne j$, the coefficient (\ref{coefficientQ}) of $q_{k,l}$ in (\ref{fij}) is unbiased, and otherwise, it is biased. This result is also the same for the coefficients of the elements of $P$ and $R$.\color{black}
\section{NUMERICAL EXPERIMENT}
\begin{figure*}
\begin{minipage}{8.5cm}
\centering
\includegraphics[width=7.5cm]{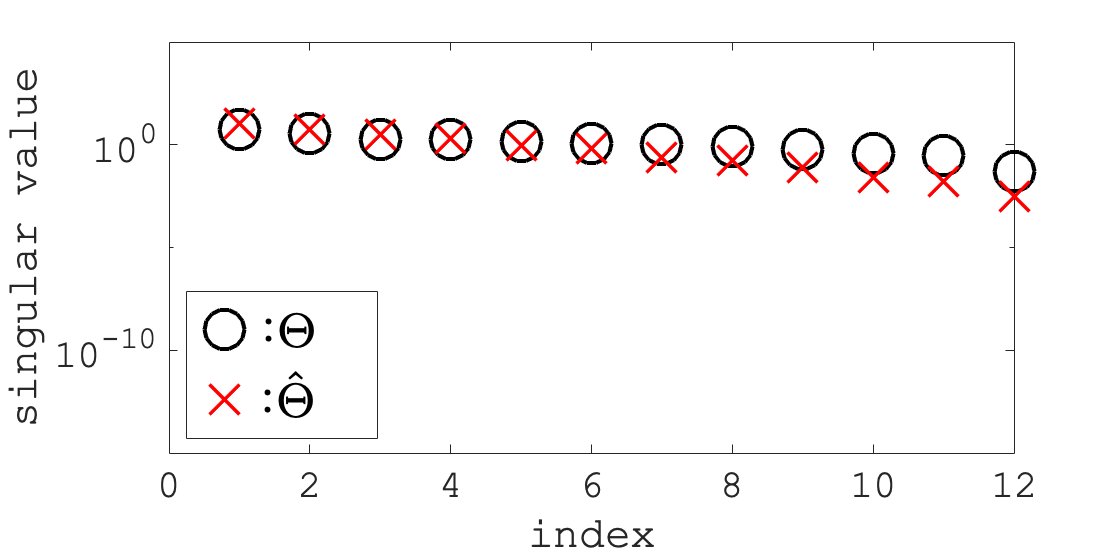}
\caption{Singular values of coefficient matrices $\Theta$ and $\hat\Theta$ in Experiment 1.}
\end{minipage}
\hfill
\begin{minipage}{8.5cm}
\centering
\includegraphics[width=7.5cm]{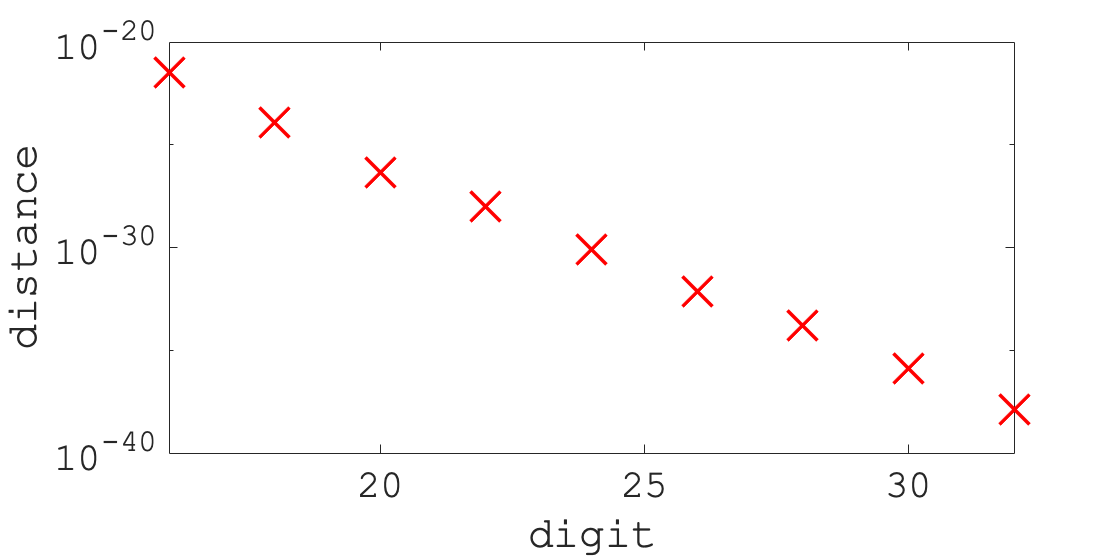}
\caption{Relationship between distance $d{\left(S,\hat S\right)}$ and number of significant digits in the computation in Experiment 1.}
\end{minipage}
\end{figure*}
\begin{figure}
\centering
\includegraphics[width=7.5cm]{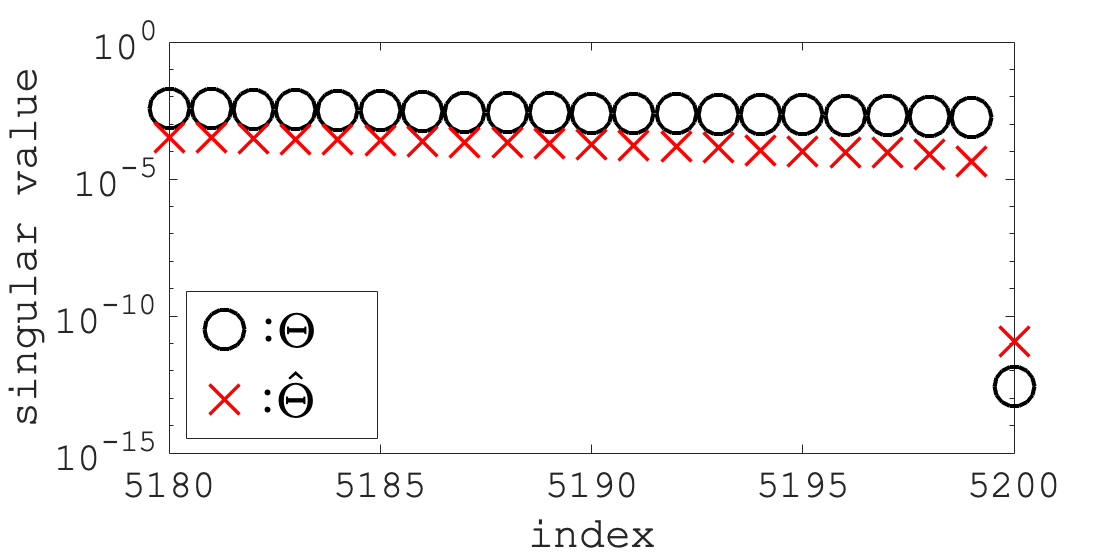}
\caption{Singular values of coefficient matrices $\Theta$ and $\hat\Theta$ in Experiment 2.}
\end{figure}
\subsection{Distance Between Solution Spaces}
In our experiments, we evaluated the ARE estimation using a distance between the solution spaces of the ARE and estimation. Let $s\in\mathbb R^{N_v}$ be a vector generated by independent elements in the matrices $P$, $Q$, and $R$ in any fixed order, where $N_v$ is defined according to the prior information. Because the ARE (\ref{ARE2}) is linear for $s$, we can transform (\ref{ARE2}) into the following equivalent form:
\begin{equation}
\label{VARE}
\Theta s=0,
\end{equation}
where $\Theta\in\mathbb R^{N_{\rm ARE}\times N_v}$ is the coefficient matrix and $N_{\rm ARE}=\frac{n\left(n+1\right)}{2}+mn$. Similarly, the estimation (\ref{estimation}) is transformed into the following equivalent form:
\begin{equation}
\label{Vestimation}
\hat\Theta s=0,
\end{equation}
where $\hat\Theta\in\mathbb R^{N_{\rm est}\times N_v}$ is the coefficient matrix and $N_{\rm est}=N_dN'_d-N'_d\left(N'_d-1\right)/2$.
We define the solution spaces $S$ and $\hat S$ of the ARE (\ref{VARE}) and estimation (\ref{Vestimation}) as follows:
\begin{equation}
S=\left\{s\in\mathbb R^{N_v}\left|\Theta s=0\right.\right\},\quad \hat S=\left\{s\in\mathbb R^{N_v}\left|\hat \Theta s=0\right.\right\}.
\end{equation}

We define the distance between the solution spaces using the approach provided \cite{book3}. Assume that $\Theta$ and $\hat\Theta$ have the same rank. Let $\Pi$, $\hat\Pi\in\mathbb R^{N_v\times N_v}$ be the orthogonal projections to $S$ and $\hat S$, respectively. The distance between $S$ and $\hat S$ is defined as follows:
\begin{equation}
\label{distance}
d{\left(S,\hat S\right)}=\left\|\Pi-\hat\Pi\right\|_2,
\end{equation}
where $\left\|\cdot\right\|_2$ is the matrix norm induced by the 2-norm of vectors. Distance $d{\left(S,\hat S\right)}$ is the maximum Euclidian distance between $\hat s\in\hat S$ and $\Pi\hat s\in S$ when $\left\|\hat s\right\|_2=1$, as explained in Fig. 2 and by the following equation:
\begin{equation}
\begin{split}
\max_{\hat s\in\hat S,\left\|\hat s\right\|_2=1}\left\|{\Pi}\hat s-\hat s\right\|_2=&\max_{\hat s\in\hat S,\left\|\hat s\right\|_2=1}\left\|\left({\Pi}-{\hat \Pi}\right)\hat s\right\|_2\\
=&d{\left(S,\hat S\right)}.
\end{split}
\end{equation}
Hence, the distance satisfies $0\le d{\left(S,\hat S\right)}\le1$.

The orthogonal projections $\Pi$ and $\hat\Pi$ are obtained from the orthogonal bases of $S$ and $\hat S$. To obtain these orthogonal bases, we use singular value decompositions of $\Theta$ and $\hat \Theta$. In our experiments, we performed the computation using MATLAB 2023a.
\subsection{Experiment 1}
In this experiment, we confirmed Theorem 2 by solving Problem 1. We set the controllable pair $\left(A,B\right)$ as follows:
\begin{equation}
A=\left[\begin{array}{*{20}{c}}
-0.2&-0.4&-0.6\\
0.4&-0.7&-0.3\\
-1.0&-0.8&-0.2
\end{array}\right],\;
B=\left[\begin{array}{*{20}{c}}
0.1&-0.6\\
-0.2&0.8\\
0.4&-0.9
\end{array}\right].
\end{equation}
The given gain $K\in\mathbb R^{2\times 3}$ is the solution of the LQR problem for the following $Q$ and $R$:
\begin{equation}
Q=\left[\begin{array}{*{20}{c}}
0.4&-0.2&0.7\\
-0.2&1.7&-0.7\\
0.7&-0.7&1.9
\end{array}\right],\;
R=\left[\begin{array}{*{20}{c}}
1.7&0.4\\
0.4&1.8
\end{array}\right].
\end{equation}
We considered the following $N_d=n+m$ observation data with $N'_d=n$ that satisfy the condition of Theorem 2:
\begin{equation}
X{\left(0\right)}=\left[\begin{array}{*{20}{c}}
0.1&0.4&0.5&-0.4&-0.1\\
0.4&0.7&1.0&0.6&0.8\\
-0.4&-1.0&0.5&-0.8&-0.4
\end{array}\right],
\end{equation}
\begin{equation}
U=\left[-K\left[x_1{\left(0\right)}\;\cdots\; x_3{\left(0\right)}\right]\begin{array}{*{20}{c}}
0.0&0.1\\
-0.9&-0.7
\end{array}\right].
\end{equation}
Using this setting, each of the ARE (\ref{ARE2}) and our estimation (\ref{estimation}) contained $12$ scalar equations for $N_v=15$ variables, that is, independent elements in the symmetric matrices $P$, $Q$, and $R$.

Fig. 2 shows the singular values of $\Theta$ and $\hat\Theta$, in descending order, that we computed using the default significant digits. As shown in Fig. 2, there is no zero singular value of $\Theta$, $\hat\Theta\in\mathbb R^{12\times 15}$. Hence, the ranks of $\Theta$ and $\hat\Theta$ are $12$, and the solution spaces $S$, $\hat S\subset\mathbb R^{15}$ are three-dimensional subspaces. To investigate the influence of the computational error, we performed the same experiment with various numbers of significant digits. Fig. 3 shows the relationship between the distance $d{\left(S,\hat S\right)}$ and the number of significant digits. As shown in Fig. 3, the logarithm of the distance decreased in proportion to the number of significant digits. Therefore, we concluded that non-zero $d{\left(S,\hat S\right)}$ was caused by the computational error and that our method correctly estimated the ARE.
\subsection{Experiment 2}
In this experiment, we confirmed that, if $Q$ and $R$ are diagonal, Theorem 1 can estimate the ARE with fewer observation data than system identification. We randomly set the matrices $A\in\mathbb R^{100\times 100}$ and $B\in\mathbb R^{100\times 50}$ by generating each element from the uniform distribution in $\left[-1,1\right]$. After the generation, we confirmed that $\left(A,B\right)$ is controllable. The given $K\in\mathbb R^{50\times 100}$ is the solution of the LQR problem for the diagonal matrices $Q$ and $R$ whose elements were generated from the uniform distribution in $\left[0.01,1\right]$. We used $N_d=N_{d\rm min}=102$ observation data that satisfy (\ref{asm_prb1}) for $N'_d=n=100$. We generated the elements of $x_{1}{\left(0\right)},\ldots,x_{102}{\left(0\right)}\in\mathbb R^{100}$ and $u_{101}$, $u_{102}\in\mathbb R^{50}$ from the uniform distribution in $\left[-1,1\right]$. The number of used data, $102$, is less than $n+m=150$ required for system identification (\ref{SI}). 

Using this setting, the ARE (\ref{ARE2}) and our estimation (\ref{estimation}) contained $N_{\rm ARE}=\frac{n\left(n+1\right)}{2}+mn=10050$ and $N_{\rm est}=N_dN'_d-N'_d\left(N'_d-1\right)/2=5250$ scalar equations, respectively. The variables of these scalar equations were $N_v=\frac{n\left(n+1\right)}{2}+n+m=5200$ independent elements in the symmetric matrix $P$ and diagonal matrices $Q$ and $R$.

Because the arbitrary-precision arithmetic was too slow to perform for the scale of this experiment, we performed the computation using the default significant digits. We indexed the singular values of $\Theta\in\mathbb R^{10050\times 5200}$ and $\hat\Theta\in\mathbb R^{5250\times 5200}$ in descending order; Fig. 4 shows the 5180--5200th singular values. From Fig. 4, we assumed that, without computational error, $\Theta$ and $\hat\Theta$ each had a zero singular value. Then, the solution spaces $S$, $\hat S\subset\mathbb R^{5200}$ were one-dimensional subspaces, and the distance $d{\left(S,\hat S\right)}$ was $4.3\times 10^{-10}$.

We compared our method with system identification. For system identification, we generated additional data $x_{103},\ldots,x_{150}$ and $u_{103},\ldots,u_{150}$ in the same manner as $x_{101}$ and $u_{101}$. Using the same approach as that for $S$, we defined the solution space $\hat S_{\rm SI}\subset\mathbb R^{N_v}$ of the ARE (\ref{ARE2}) using $A$, $B$, and $K$ identified by (\ref{SI}). Then, the distance $d{\left(S,\hat S_{\rm SI}\right)}$ was $2.8\times 10^{-10}$. Therefore, our method estimated the ARE with the same order of error as system identification using fewer observation data than system identification.
\subsection{Experiment 3}
In this experiment, we compared our method and system identification under a practical conditions: a single noisy and sequential data and sparse but not diagonal $Q$ and $R$. We set the matrices $A\in\mathbb R^{40\times 40}$ and $B\in\mathbb R^{40\times 20}$ in the same manner as Experiment 2. The control gain $K\in\mathbb R^{20\times 40}$ is the solution of the LQR problem for the sparse matrices $Q$ and $R$. We performed three operations to generate $Q$ and $R$. First, we generated matrices $M_Q\in\mathbb R^{40\times 40}$ and $M_R\in\mathbb R^{20\times 20}$ in the same manner as $A$ and set $Q:=M_Q^\top M_Q$ and $R:=M_R^\top M_R$. Second, we randomly set $800$ non-diagonal elements of $Q$ to $0$ to make it sparse while maintaining symmetry. Similarly, we set $200$ elements of $R$ to $0$. Finally, add the product of a constant and the identity matrix to $Q$ and $R$ so that the maximum eigenvalues are $10$ times as large as the minimum ones, respectively.

\begin{figure*}
\begin{minipage}{8.5cm}
\centering
\includegraphics[width=7.5cm]{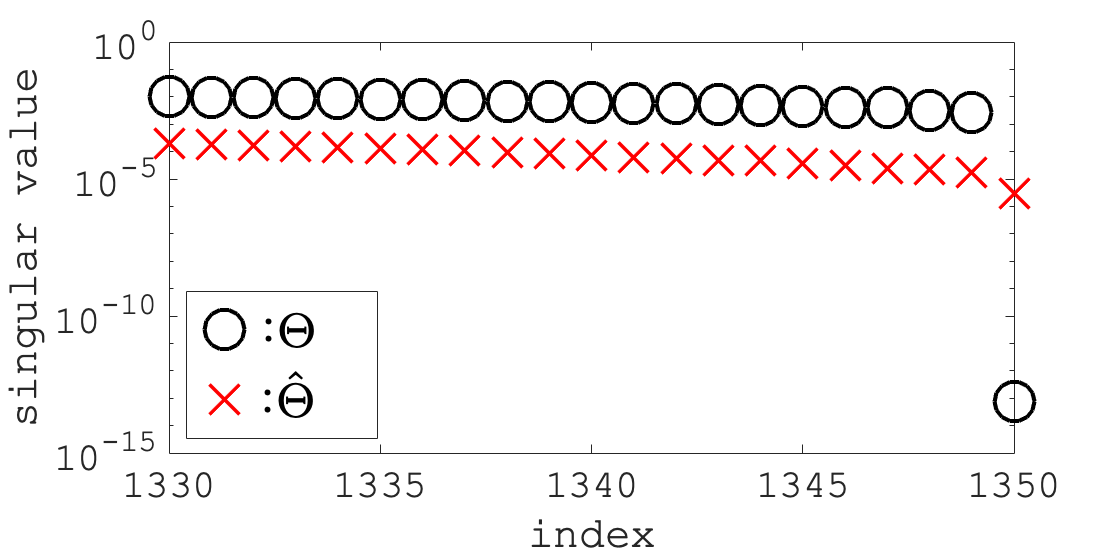}
\caption{Singular values of coefficient matrices $\Theta$ and $\hat\Theta$ in Experiment 3.}
\end{minipage}
\hfill
\begin{minipage}{8.5cm}
\centering
\includegraphics[width=7.5cm]{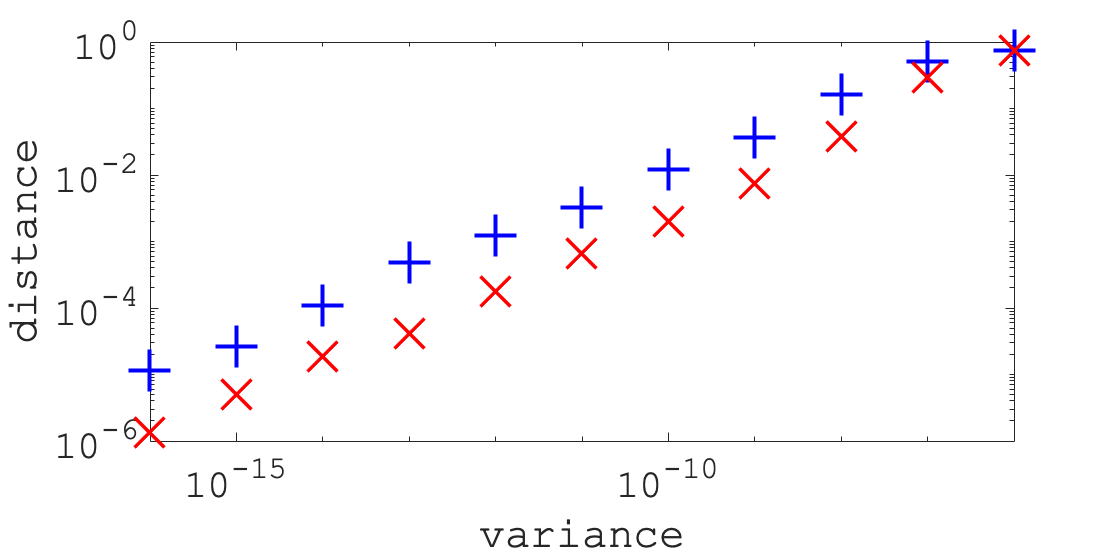}
\caption{Relationship between the noise variance $\sigma^2$ and distances $d{\left(S,\hat S\right)}$ (red ``$\times$'') and $d{\left(S,\hat S_{\rm SI}\right)}$ (blue ``$+$'') in Experiment 3.}
\end{minipage}
\end{figure*}

To generate data, we used the following system:
\begin{equation}
\label{noisy_system}
\begin{split}
x^*{\left(k+1\right)}=&Ax^*{\left(k\right)}+Bu^*{\left(k\right)},\\
u^*{\left(k\right)}=&-Kx^*{\left(k\right)}+v{\left(k\right)},
\end{split}
\end{equation}
where $x^*{\left(k\right)}\in{\mathbb R}^n$ and $u^*{\left(k\right)}\in{\mathbb R}^m$ are the state and input at time $k\in\mathbb Z_+$, respectively. If $\left\|x^*{\left(k\right)}\right\|_2\le1$, the vector $v{\left(k\right)}\in{\mathbb R}^n$ with norm $0.2$ is the product of a scalar constant and a random vector whose elements follow a uniform distribution in $\left[-1,1\right]$. Otherwise, $v{\left(k\right)}=0$. We ran the system (\ref{noisy_system}) through $N_d=200$ time steps with a random initial state $x^*{\left(0\right)}$ satisfying $\left\|x^*{\left(0\right)}\right\|_2\le1$ and obtained data as follows:
\begin{equation}
\begin{split}
x_k{\left(0\right)}=&x^*{\left(k-1\right)}+\varepsilon_{\rm state}{\left(k-1\right)},\\
u_k=&u^*{\left(k-1\right)}+\varepsilon_{\rm input}{\left(k-1\right)},\\
x_k{\left(1\right)}=&x^*{\left(k\right)}+\varepsilon_{\rm state}{\left(k\right)},
\end{split}
\end{equation}
where $\varepsilon_{\rm state}\in{\mathbb R}^n$ and $\varepsilon_{\rm input}\in{\mathbb R}^m$ are the observation noise whose elements follow a normal distribution with a mean $0$ and variance $\sigma^2$. We explain the value of $\sigma^2$ later. Throughout the simulation, $v{\left(k\right)}=0$ holds $121$ times. Hence, $N'_d=121$. We sorted the data $\left(x_i{\left(0\right)},u_i,x_i{\left(1\right)}\right)$ $\left(i\in\left\{1,\ldots,N_d\right\}\right)$ so that (\ref{asm_prb1}) holds if noise does not exist.

Under the above condition, the ARE (\ref{ARE2}) and our estimation (\ref{estimation}) contained $N_{\rm ARE}=\frac{n\left(n+1\right)}{2}+mn=1620$ and $N_{\rm est}=N_dN'_d-N'_d\left(N'_d-1\right)/2=16940$ scalar equations, respectively. The variables of these scalar equations were $N_v=n\left(n+1\right)+\frac{m\left(m+1\right)}{2}-500=1350$ independent elements in the symmetric matrix $P$ and sparse matrices $Q$ and $R$.

We indexed the singular values of $\Theta\in\mathbb R^{1620\times 1350}$ and $\hat\Theta\in\mathbb R^{16940\times 1350}$ in descending order. Fig. 5 shows the 1330--1350th singular values when $\sigma^2=10^{-8}$. Although there is no zero singular value due to noise, the last singular value is noticeably small as shown in Fig. 5. Hence, we can conclude that the solution spaces $S$, $\hat S\subset\mathbb R^{1350}$ were one-dimensional subspaces.

We conducted experiments with different noise variance $\sigma^2$ from $10^{-6}$ to $10^{-16}$. Also, the noise seed differs for each experiment, but the other seeds are the same. Fig. 6 shows the relationship between the noise variance $\sigma^2$ and distances $d{\left(S,\hat S\right)}$ and $d{\left(S,\hat S_{\rm SI}\right)}$. As shown in Fig. 6, our method outperformed system identification in almost all experiments. In the experiments of variances from $10^{-8}$ to $10^{-16}$, the ratio of $d{\left(S,\hat S\right)}$ to $d{\left(S,\hat S_{\rm SI}\right)}$ is approximately constant, and the average ratio is $0.17$. Because the maximum value of the distances is $1$, we conclude that the estimations using both methods failed with a larger noise variance than $10^{-8}$.\color{black}
\section{CONCLUSIONS}
In this paper, we proposed a method to estimate the ARE with respect to an unknown discrete-time system from the input and state observation data. Our method transforms the ARE into a form calculated without the system model by multiplying the observation data on both sides. We proved that our estimated equation is equivalent to the ARE if the data are sufficient for system identification. The main feature of our method is the direct estimation of the ARE without identifying the system. This feature enables us to economize the observation data using prior information about the objective function. We conducted a numerical experiment that demonstrated that that our method requires less data than system identification if the prior information is sufficient. 
\addtolength{\textheight}{-3cm}

\newpage\onecolumn\setcounter{equation}{0}
\section*{SUPPLEMENTARY MATERIAL}
We give another proof of Theorem 1 using the following Bellman equation:
\begin{equation}
\label{bellmanSM}
\begin{split}
x{\left(k\right)}^\top Px{\left(k\right)}=&x{\left(k+1\right)}^\top Px{\left(k+1\right)}+x{\left(k\right)}^\top Qx{\left(k\right)}\\
&+x{\left(k\right)}^\top K^\top RKx{\left(k\right)}.
\end{split}
\end{equation}
To this end, we prove the following theorem equivalent to Theorem 1.

{\bf Theorem 1': }{\it Consider Problem 1. Let $i$, $j\in\left\{1,\ldots,N_d\right\}$ be arbitrary natural numbers. Suppose that the following ARE holds:
\begin{equation}
\begin{split}
A^\top PA-P+Q-K^\top\left(R+B^\top PB\right)K=&0,\\
B^\top PA-\left(R+B^\top PB\right)K=&0.
\end{split}
\end{equation}
Then, we have $f_{i,j}=0$ if $u_i=-Kx_i{\left(0\right)}$ holds, where $f_{i,j}$ is defined as
\begin{equation}
\label{fijSM}
f_{i,j}=x_i{\left(1\right)}^\top Px_j{\left(1\right)}+x_i{\left(0\right)}^\top \left(Q-P\right)x_j{\left(0\right)}+u_i^\top Ru_j.
\end{equation}
}

\begin{proof}
From the Bellman equation (\ref{bellmanSM}), we can readily obtain $f_{i,j}=0$ if $i=j$. Next, suppose $i \neq j$ and that $u_j=-Kx_j{\left(0\right)}$ also holds. Then, we have $\left[x_i{\left(1\right)}+x_j{\left(1\right)}\right]=A\left[x_i{\left(0\right)}+x_j{\left(0\right)}\right]+B\left(u_i+u_j\right)$ and $u_i+u_j=-K\left[x_i{\left(0\right)}+x_j{\left(0\right)}\right]$. Hence, from the Bellman equation (\ref{bellmanSM}), we obtain
\begin{equation}
\begin{split}
0=&\left[x_i{\left(1\right)}+x_j{\left(1\right)}\right]^\top P\left[x_i{\left(1\right)}+x_j{\left(1\right)}\right]\\
&+\left[x_i{\left(0\right)}+x_j{\left(0\right)}\right]^\top\left(Q-P\right)\left[x_i{\left(0\right)}+x_j{\left(0\right)}\right]\\
&+\left(u_i+u_j\right)^\top R\left(u_i+u_j\right)\\
&=f_{i,i}+f_{i,j}+f_{j,i}+f_{j,j}.
\end{split}
\end{equation}
\samepage Because $f_{i,i}=f_{j,j}=0$ and $f_{i,j}=f_{j,i}$, we obtain
\begin{equation}
\label{231024133217}
\left[u_i=-Kx_i{\left(0\right)},u_j=-Kx_j{\left(0\right)}\right]\Rightarrow f_{i,j}=0.
\end{equation}

Finally, we consider the case without the assumption $u_j=-Kx_j{\left(0\right)}$. Let $\Delta u_j=u_j+Kx_j{\left(0\right)}$. Substituting $x_j{\left(1\right)}=Ax_j{\left(0\right)}+Bu_j$ and $u_j=\Delta u_j-Kx_j{\left(0\right)}$ to (\ref{fijSM}), we have
\begin{equation}
\label{231024143626}
\begin{split}
f_{i,j}=&x_i{\left(1\right)}^\top P\left(A-BK\right)x_j{\left(0\right)}+x_i{\left(1\right)}^\top PB\Delta u_j\\
&+x_i{\left(0\right)}^\top\left(Q-P\right)x_j{\left(0\right)}-u_i^\top RKx_j{\left(0\right)}+u_i^\top R\Delta u_j.
\end{split}
\end{equation}
From (\ref{231024133217}), the sum of the first, third, and fourth terms on the right-hand side of (\ref{231024143626}) is zero. Hence, substituting $x_i{\left(1\right)}=Ax_i{\left(0\right)}+Bu_i$ and $u_i=-Kx_i{\left(0\right)}$ to (\ref{231024143626}), we have
\begin{equation}
f_{i,j}=x_i{\left(0\right)}^\top\left[-K^\top\left(B^\top PB+R\right)+A^\top PB\right]\Delta u_j.
\end{equation}
From $K=\left(R+B^\top PB\right)^{-1}B^\top PA$, we obtain $f_{i,j}=0$, which proves Theorem 1'.
\end{proof}

\begin{thebibliography}{99}
\bibitem{Mombaur}
K. Mombaur, A. Truong, and J.-P. Laumond, ``From human to humanoid locomotion--an inverse optimal control approach'', {\it Autonomous Robots}, vol. 28, pp. 369--383, 2010.

\bibitem{Li}
W. Li, E. Todorov, and D. Liu, ``Inverse Optimality Design for Biological Movement Systems'', {\it IFAC Proceedings Volumes}, vol. 44, no. 1, pp. 9662--9667, 2011.

\bibitem{El-Hussieny}
H. El-Hussieny, A.A. Abouelsoud, S.F.M. Assal, and S.M. Megahed, ``Adaptive learning of human motor behaviors: An evolving inverse optimal control approach'', {\it Engineering Applications of Artificial Intelligence}, vol. 50, pp. 115--124, 2016

\bibitem{Alexander}
R.M. Alexander, ``The gaits of bipedal and quadrupedal animals'', {\it International Journal of Robotics Research}, 3(2), pp. 49--59, 1984.

\bibitem{Kalman}
R.E. Kalman, ``When Is a Linear Control System Optimal?'', {\it Journal of Basic Engineering}, 86(1), pp. 51--60, 1964.

\bibitem{Anderson}
B.D.O. Anderson, {\it The Inverse Problem of Optimal Control, Technical report: Stanford University, vol. 6560, no.3}, Stanford University, 1966.

\bibitem{Molinari}
B. Molinari, ``The stable regulator problem and its inverse'', {\it IEEE Transactions on Automatic Control}, vol. 18, no. 5, pp. 454--459, 1973.

\bibitem{Moylan}
P. Moylan and B. Anderson, ``Nonlinear regulator theory and an inverse optimal control problem'', {\it IEEE Transactions on Automatic Control}, vol. 18, no. 5, pp. 460--465, 1973.

\bibitem{Ng}
A.Y. Ng and S. Russell, ``Algorithms for Inverse Reinforcement Learning'', {\it Proceedings of the Seventeenth International Conference on Machine Learning}, pp. 663--670, 2000.

\bibitem{Ab Azar}
N. Ab Azar, A. Shahmansoorian, and M. Davoudi, ``From inverse optimal control to inverse reinforcement learning: A historical review'', {\it Annual Reviews in Control}, vol. 50, pp. 119--138, 2020.

\bibitem{Priess}
M.C. Priess, R. Conway, J. Choi, J.M. Popovich, and C. Radcliffe, ``Solutions to the Inverse LQR Problem With Application to Biological Systems Analysis'', {\it IEEE Transaction on Control Systems Technology}, vol. 23, no. 2, pp. 770--777, 2015.


\bibitem{Herman}
M. Herman, T. Gindele, J. Wagner, F. Schmitt, and W. Burgard, ``Inverse reinforcement learning with simultaneous estimation of rewards and dynamics''. {\it Artificial intelligence and statistics}, PMLR, pp. 102--110, 2016.

\bibitem{Aghasadeghi}
N. Aghasadeghi and T. Bretl, ``Maximum entropy inverse reinforcement learning in continuous state spaces with path integrals''. {\it IEEE/RSJ International Conference on Intelligent Robots and Systems}, pp. 1561--1566, 2011.

\bibitem{continuous}
S. Sugiura, R. Ariizumi, M. Tanemura, T. Asai, and S. Azuma, {\it Data-driven Estimation of Algebraic Riccati Equation for Inverse Linear Quadratic Regulator Problem}, SICE Annual Conference, 2023.

\bibitem{book1}
P.J. Antsaklis and A.N. Michel, {\it Linear Systems}, Birkh\"auser Boston, 2005.

\bibitem{book2}
M. bardi, {\it Optimal Control and Viscosity Solutions of Hamilton-Jacobi-Bellman Equations}, Birkh\"auser, 1997.

\bibitem{book3}
G.H. Golub, and C.F. Van Loan, {\it Matrix Computations}, JHU Press, 1997.
\end{thebibliography}
\end{document}